\documentclass[a4paper,11pt,twoside]{article}
\RequirePackage[a4paper,head=1cm]{geometry}
\usepackage[utf8x]{inputenc}
\usepackage{textcomp}
\usepackage{amsmath,amsfonts,verbatim,afterpage,euscript,mathrsfs,amssymb,empheq}
\usepackage{amsthm}
\usepackage{dsfont}
\usepackage{hyperref}
\usepackage{pst-fill,pst-grad}
\usepackage{textcomp}
\usepackage{multicol}
\usepackage{amsmath}
\usepackage{amssymb}
\usepackage{hyperref}
\usepackage{dsfont}

\newtheorem{Theoreme}{Theorem}
\newtheorem{Lemme}{Lemma}[section]
\newtheorem{Corollaire}{Corollary}[section]

\newtheorem{Remarque}{\bf Remark}
\newcommand{\mysection}{\setcounter{equation}{0} \section}

\title{A new pointwise inequality for rough operators and applications} 
\author{Diego Chamorro\footnote{Laboratoire de Math\'ematiques et Mod\'elisation d'Evry (LaMME) - UMR 8071. Universit\'e d'Evry Val d'Essonne, 23 Boulevard de France, 91037 Evry Cedex, France. email: \textit{diego.chamorro@univ-evry.fr}}, Anca-Nicoleta Marcoci\footnote{Department of Mathematics and Computer Science. Technical University of Civil Engineering, Bucharest, Bld. Lacul Tei, no. 124, sector 2. Romania. email: \textit{anca.marcoci@utcb.ro}}, Liviu-Gabriel Marcoci\footnote{Department of Mathematics and Computer Science. Technical University of Civil Engineering, Bucharest, Bld. Lacul Tei, no. 124, sector 2. Romania.  email: \textit{liviu.marcoci@utcb.ro}}.} 
\begin{document} 
\maketitle 
\begin{scriptsize}
\abstract{ \noindent We study in this article a new pointwise estimate for ``rough'' singular integral operators. From this pointwise estimate we will derive Sobolev type inequalities in a variety of functional spaces.}\\

\end{scriptsize}

\mysection{Introduction}
When dealing with operators and functions, pointwise estimates (in addition to giving a deep understanding of the properties of the operators considered) usually provide many interesting inequalities. A particular example that will guide our study is the following: for $n\geq 2$, if $f:\mathbb{R}^n\longrightarrow \mathbb{R}$ is a function such that $f\in \mathcal{C}^\infty_0(\mathbb{R}^n)$, then we have the estimate
\begin{equation}\label{Ineq_IntroRiesz}
|f(x)|\leq CI_{1}(|\nabla f|)(x),
\end{equation}
 (see \cite[Lemma 3.3]{Kinnunen} for a proof) where the operator $I_1$ corresponds to the usual Riesz potential defined by the expression
$$I_1(f)(x)= \pi ^{-n/2}2^{-1}\frac{\Gamma ((n-1 )/2)}{\Gamma (1/2)}\int_{\mathbb{R}^n}\frac{f(y)}{|x-y|^{n-1}}dy.$$
Note now that since the Riesz potentials $I_1$ are bounded from $L^p(\mathbb{R}^n)$ to $L^{q}(\mathbb{R}^n)$ with $1<p<n$ and $\frac{1}{q}=\frac{1}{p}-\frac{1}{n}$ (see \cite[Theorem 6.1.3]{Grafakos}), we can easily deduce from the previous pointwise estimate (\ref{Ineq_IntroRiesz}) the functional inequality
$$\|f\|_{L^{q}}\leq C\|I_{1}(|\nabla f|)\|_{L^{q}}\leq C\|\nabla f\|_{L^p},$$
which is nothing but the classical Sobolev inequality $\|f\|_{L^{q}}\leq C\|\nabla f\|_{L^p}$.\\

Of course, the pointwise estimate (\ref{Ineq_IntroRiesz}) as well as the previous Sobolev inequality admit several modifications (and different proofs) and in the recent works \cite{Hoang1}, \cite{Hoang} and \cite{Li} the following generalization of the inequality (\ref{Ineq_IntroRiesz}) was studied
\begin{equation}\label{Ineq_IntroOper}
|T(f)(x)|\leq CI_{1}(|\nabla f|)(x),
\end{equation}
where the operators $T$ considered in these references are iterates of the Hardy-Littlewood maximal function, spherical maximal operator, $L^r$-maximal operators, Lorentz-based maximal operators and ``rough'' singular integral operators. \\

In this article we are going to focus our study on the pointwise estimate (\ref{Ineq_IntroOper}) with \emph{rough singular integral operators} which are defined as follows: for a locally integrable function $f:\mathbb{R}^n\longrightarrow \mathbb{R}$, we will consider the operator $T_\Omega$ associated to a function $\Omega:\mathbb{S}^{n-1}\longrightarrow \mathbb{R}$ by the expression
\begin{equation}\label{Def_Operator}
T_\Omega(f)(x)=p.v.\int_{\mathbb{R}^n}\frac{\Omega(y/ |y|)}{|y|^n}f(x-y)dy.
\end{equation}
The properties of the function $\Omega$ are absolutely essential to understand the behavior of the associated operator $T_\Omega$ in connection to the estimate (\ref{Ineq_IntroOper}). Indeed, by considering a function $\Omega$ such that $\Omega \in L^1(\mathbb{S}^{n-1})$, $\displaystyle{\int_{\mathbb{S}^{n-1}}\Omega \ d\sigma=0}$ and such that $\Omega \in L^{\infty}(\mathbb{S}^{n-1})$ then the following estimate was proven in \cite{Li}:
\begin{equation}\label{Ineq_IntroOper1}
|T_{\Omega}(f)(x)|\leq C\|\Omega\|_{L^{\infty}(\mathbb{S}^{n-1})}I_{1}(|\nabla f|)(x)\qquad \mbox{for }  f\in\mathcal{C}^\infty_0(\mathbb{R}^n).
\end{equation}
This pointwise estimate is of particular interest when considering Sobolev inequalities: by the boundedness properties of the Riesz potential $I_1$ we can easily deduce the inequality 
\begin{equation}\label{Ineq_IntroSobolevOper0}
\|T_{\Omega}(f)\|_{L^{q}}\leq C\|\Omega\|_{L^{\infty}(\mathbb{S}^{n-1})}\|I_{1}(|\nabla f|)\|_{L^{q}}\leq C\|\Omega\|_{L^{\infty}(\mathbb{S}^{n-1})}\|\nabla f\|_{L^{p}},
\end{equation}
with $1<p<n$ and $\frac{1}{q}=\frac{1}{p}-\frac{1}{n}$. Estimates of the form (\ref{Ineq_IntroOper1}) also allow several weighted versions of the previous Sobolev-like inequalities:  for example if the Riesz potential $I_1$ is bounded from $L^p(u)$ to $L^q(v)$ where $1<p,q<+\infty$ and where $u,v$ are suitable weights, then it is possible to derive from (\ref{Ineq_IntroOper1}) the inequality
\begin{equation}\label{Ineq_IntroSobolevOper1}
\|T_{\Omega}(f)\|_{L^q(v)}\leq C\|\Omega\|_{L^{\infty}(\mathbb{S}^{n-1})}\|I_{1}\|_{\mathcal{B}(L^p(u), L^q(v))}\|\nabla f\|_{L^p(u)}.
\end{equation}
More recently, in \cite{Hoang}, the estimate (\ref{Ineq_IntroOper1}) (and consequently the inequalities (\ref{Ineq_IntroSobolevOper0}) or (\ref{Ineq_IntroSobolevOper1})) was improved by considering the more general condition $\Omega \in L^{n,\infty}(\mathbb{S}^{n-1})$ where the space $L^{n,\infty}$ is a Lorentz space (recall that since $\sigma(\mathbb{S}^{n-1})<+\infty$ we have $L^\infty(\mathbb{S}^{n-1})\subset L^n(\mathbb{S}^{n-1})$ and that we always have the space inclusion $L^n(\mathbb{S}^{n-1})\subset L^{n,\infty}(\mathbb{S}^{n-1})$). The pointwise control is then the following:
\begin{equation}\label{Ineq_IntroOperLorentz}
|T_{\Omega}(f)(x)|\leq C\|\Omega\|_{L^{n,\infty}(\mathbb{S}^{n-1})}I_{1}(|\nabla f|)(x),
\end{equation}
from which we deduce the inequalities
\begin{eqnarray}
\|T_{\Omega}(f)\|_{L^q}&\leq &C\|\Omega\|_{L^{n,\infty}(\mathbb{S}^{n-1})}\|\nabla f\|_{L^p}\qquad \mbox{with } 1<p<n \mbox{ and }  \tfrac{1}{q}=\tfrac{1}{p}-\tfrac{1}{n}\quad \mbox{ or }  \label{Ineq_IntroSobolevOper2}\\[3mm]
\|T_{\Omega}(f)\|_{L^q(v)}&\leq &C\|\Omega\|_{L^{n,\infty}(\mathbb{S}^{n-1})}\|I_{1}\|_{\mathcal{B}(L^p(u), L^q(v))}\|\nabla f\|_{L^p(u)},\notag
\end{eqnarray}
for suitable weights $u,v$ and adapted weighted spaces $L^q(v)$ and $L^p(u)$. 
Let us mention that weak endpoints were also considered in the references \cite{Li} and \cite{Hoang} as well as many other consequences of the pointwise inequalities (\ref{Ineq_IntroOper1}) and (\ref{Ineq_IntroOperLorentz}).\\

Remark now that in \cite{Hoang1}, an interesting modification of the pointwise estimates (\ref{Ineq_IntroOper1}) and (\ref{Ineq_IntroOperLorentz}) was studied. Indeed, for $0<\mathfrak{p}<n$ we can consider the Riesz potential $I_\mathfrak{p}$ defined by 
\begin{equation}\label{Ineq_IntroRieszAlpha}
I_\mathfrak{p}(f)(x)= \pi ^{-n/2}2^{-\mathfrak{p}}\frac{\Gamma ((n-\mathfrak{p} )/2)}{\Gamma (\mathfrak{p}/2)}\int_{\mathbb{R}^n}\frac{f(y)}{|x-y|^{n-\mathfrak{p}}}dy,
\end{equation}
and we can define the operator $T_{\Omega, \mathfrak{p}}$ by 
$$T_{\Omega,\mathfrak{p}}(f)(x)=p.v.\int_{\mathbb{R}^n}\frac{\Omega(y/ |y|)}{|y|^{n+1-\mathfrak{p}}}f(x-y)dy.$$
Thus, if $\Omega \in L^1(\mathbb{S}^{n-1})$, $\displaystyle{\int_{\mathbb{S}^{n-1}}\Omega\ d\sigma=0}$ and $\Omega \in L^{n,\infty}(\mathbb{S}^{n-1})$, it was proven in \cite[Theorem 1.1]{Hoang1} the following pointwise estimate
$$|T_{\Omega,\mathfrak{p}}(f)(x)|\leq C\|\Omega\|_{L^{n,\infty}(\mathbb{S}^{n-1})}I_{\mathfrak{p}}(|\nabla f|)(x), \qquad \mbox{for } f\in\mathcal{C}^\infty_0(\mathbb{R}^n),$$
from which several functional inequalities of Sobolev-type of the type (\ref{Ineq_IntroSobolevOper0}) or (\ref{Ineq_IntroSobolevOper1}) are deduced. \\

In this work we will start by proving a variation of the pointwise estimate (\ref{Ineq_IntroOperLorentz}) where we will introduce two modifications. First we will replace the boundedness information of the function $\Omega$ stated in terms of the Lorentz space $L^{n,\infty}(\mathbb{S}^{n-1})$ by a more general condition given by $\Omega\in L^{\rho}(\mathbb{S}^{n-1})$ with $1<\rho<n$. Indeed, in this case and since $\sigma(\mathbb{S}^{n-1})<+\infty$, we have $L^{n,\infty}(\mathbb{S}^{n-1})\subset L^{\rho}(\mathbb{S}^{n-1})$. The second modification consists in replacing the Riesz potential $I_1$ in (\ref{Ineq_IntroOperLorentz}) by a mixed information which involves the Hardy-Littlewood maximal function $\mathscr{M}$ defined by 
$$\mathscr{M}(f)(x)=\displaystyle{\underset{B \ni x}{\sup } \;\frac{1}{|B|}\int_{B }|f(y)|dy},$$
and a Morrey space $\dot{\mathcal{M}}^{p,q}(\mathbb{R}^n)$ defined for $1\leq p\leq q<+\infty$ by the condition
$$\|f\|_{\dot{\mathcal{M}}^{p,q}}=\underset{x\in \mathbb{R}^n, \; r>0}{\sup}\left(\frac{1}{r^{n(1-\frac{p}{q})}}\int_{B(x,r)}|f(y)|^p dy\right)^{\frac{1}{p}}<+\infty,$$
see Section \ref{Secc_FuncSpaces} below for more details on Morrey spaces. As we shall see, these modifications of the pointwise estimate (\ref{Ineq_IntroOperLorentz}) will provide an interesting framework from which we will deduce new functional inequalities.\\

In this context, our main result reads as follows: 
\begin{Theoreme}[{\bf Pointwise inequality}]\label{Theo_Pointwise}
Over the space $\mathbb{R}^n$ with $n\geq 2$, consider $\Omega$ a function such that $\Omega \in L^1(\mathbb{S}^{n-1})$, $\displaystyle{\int_{\mathbb{S}^{n-1}}\Omega \ d\sigma=0}$ and such that $\Omega\in L^\rho(\mathbb{S}^{n-1})$ with $1<\rho<n$ and consider the operator $T_{\Omega}$ associated to the function $\Omega$ as defined in (\ref{Def_Operator}). \\

\noindent If $f:\mathbb{R}^n\longrightarrow \mathbb{R}$ is a regular function such that $\nabla f\in \dot{\mathcal{M}}^{\mathfrak{p},\mathfrak{q}}(\mathbb{R}^n)$ where $1< \mathfrak{p}\leq \mathfrak{q}<n$ and where $\mathfrak{p}$ is a real parameter such that $1<\overline{\rho}=\frac{\rho n}{\rho n+\rho-n}\leq \mathfrak{p}<n$, then we have the following pointwise estimate
\begin{equation}\label{pointwise_estimate}
|T_\Omega(f)(x)|\leq C \|\Omega\|_{L^\rho(\mathbb{S}^{n-1})}\big(\mathscr{M}(|\nabla f|^{\mathfrak{p}})(x)\big)^{\frac{1}{\mathfrak{p}}-\frac{\mathfrak{q}}{\mathfrak{p}n}}\|\nabla f\|_{\dot{\mathcal{M}}^{\mathfrak{p},\mathfrak{q}}}^{\frac{\mathfrak{q}}{n}}.
\end{equation}
\end{Theoreme}
\noindent Some remarks are in order. Note that since our first motivation was to work with more general operators and to extend the condition $\Omega\in L^{n,\infty}(\mathbb{S}^{n-1})$ in the estimate (\ref{Ineq_IntroOperLorentz}), it was quite natural to consider the  Lebesgue spaces $L^{\rho}(\mathbb{S}^{n-1})$ and thus in order to obtain the space inclusion $L^{n,\infty}(\mathbb{S}^{n-1})\subset L^{\rho}(\mathbb{S}^{n-1})$, we require the condition $\rho<n$ (see formula (\ref{InclusionLorentzLpUloc}) below). Note in particular that the value $\rho=n$ is allowed (and then we can set $\mathfrak{p}=1$) but in this case the corresponding space for the function $\Omega$ would be $L^{n}(\mathbb{S}^{n-1})$, which would be an improvement of (\ref{Ineq_IntroOper1}) but not of (\ref{Ineq_IntroOperLorentz}) and for this reason we will restrict ourselves to the case $1<\rho<n$.  To continue, let us remark now that the boundedness of the Hardy-Littlewood maximal function $\mathscr{M}$ is in many situations  (\emph{i.e.} outside of usual cases) easier to establish than the boundedness of the Riesz potentials: indeed, in the case of the maximal function the same functional space can be naturally considered (\emph{i.e.} we have $\|\mathscr{M}(f)\|_{E}\leq C\|f\|_E$) while the boundedness of Riesz potentials always involves two spaces (\emph{i.e.} we have $\|I_{1}(f)\|_{F}\leq C \|I_{1}\|_{\mathcal{B}(E, F)} \|f\|_E$) and this fact will make the estimate (\ref{pointwise_estimate}) more robust than (\ref{Ineq_IntroOperLorentz}). We will also see in the lines below that the pointwise estimate (\ref{pointwise_estimate}) will provide sharper Sobolev-type estimates than (\ref{Ineq_IntroOperLorentz}). Let us mention finally that we can actually obtain a slightly sharper estimate than (\ref{pointwise_estimate}) by considering the operator $T^*_{\Omega}$ defined in the expression (\ref{Def_OperatorTstar}). Indeed, the control (\ref{pointwise_estimate}) will be a consequence of the inequality (\ref{EstimationInitialRefined1}) below. However, and for the sake of simplicity, we do not introduce this operator $T^*_{\Omega}$ in the statement of our Theorem \ref{Theo_Pointwise}.\\

We explore now very natural applications of the pointwise estimate (\ref{pointwise_estimate}). Indeed, we have:
\begin{Theoreme}[{\bf Refined Sobolev inequalities}]\label{Theo_RefinedSobolevInequalities}
Over the space $\mathbb{R}^n$ with $n\geq 2$, consider $\Omega$ a function such that $\Omega \in L^1(\mathbb{S}^{n-1})$, $\displaystyle{\int_{\mathbb{S}^{n-1}}\Omega \ d\sigma=0}$ and such that $\Omega\in L^\rho(\mathbb{S}^{n-1})$ with $1<\rho<n$ and consider the operator $T_{\Omega}$ associated to the function $\Omega$ as defined in (\ref{Def_Operator}). \\

\noindent If $f:\mathbb{R}^n\longrightarrow \mathbb{R}$ is a regular function such that $\nabla f\in \dot{\mathcal{M}}^{\mathfrak{p},\mathfrak{q}}(\mathbb{R}^n)$ where $1< \mathfrak{p}\leq \mathfrak{q}<n$ and where $\mathfrak{p}$ is a real parameter such that $1<\overline{\rho}=\frac{\rho n}{\rho n+\rho-n} \leq  \mathfrak{p}$ and if we have $\nabla f\in  L^{\sigma}(\mathbb{R}^n)$ for some index $\sigma$ such that $\mathfrak{p}<\sigma<+\infty$, then for an index $r=\frac{n\sigma}{n-\mathfrak{q}}>1$ and a power $\theta=\frac{\mathfrak{q}}{n}$, we have the inequality
\begin{equation}\label{InequalitySobolev_Intro}
\|T_\Omega(f)\|_{L^r}\leq C_\Omega \left\|\nabla f\right\|_{L^{\sigma}}^{1-\theta}\|\nabla f\|_{\dot{\mathcal{M}}^{\mathfrak{p},\mathfrak{q}}}^{\theta}.
\end{equation}
\end{Theoreme}

\noindent These theorems constitute the core of our article, but in the sections below we will also study weighted inequalities and some functional estimates in different frameworks.\\

To conclude this introduction, we point out that in this work we do not consider weak endpoints of the inequalities of the type (\ref{InequalitySobolev_Intro}) as this will require a different treatment. \\

The plan of the article is the following. In Section \ref{Secc_FuncSpaces} we will recall the definitions and the main properties of the functional spaces used here. In Section \ref{Secc_ProofTh1} we will prove Theorem \ref{Theo_Pointwise} and Section \ref{Secc_ProofTheo2} will be devoted to the proof of Theorem \ref{Theo_RefinedSobolevInequalities}. In Section \ref{Secc_Weights} we shall also present some weighted variants of the inequality (\ref{InequalitySobolev_Intro}) while in Section \ref{Secc_Orlicz} we will extend the inequality (\ref{InequalitySobolev_Intro}) to the framework of Orlicz spaces. Finally, in Section \ref{Secc_ClassicalLorentz} we will consider the framework of classical Lorentz spaces. 
\mysection{Some functional spaces and classical inequalities}\label{Secc_FuncSpaces}
In this section we recall the definitions and some well known properties of the functional spaces that will be used here. 
\begin{itemize}
\item For $1\leq p<+\infty$ and for $A=\mathbb{R}^n$ or $A \subset \mathbb{R}^n$, the usual Lebesgue space $L^p(A)$ are defined by the classical condition
$\|f\|_{L^p}=\displaystyle{\left(\int_{A}|f(x)|^pdx\right)^{\frac{1}{p}}}<+\infty$. Recall in particular that if $A$ is a bounded subset then we have the space inclusions $L^{p_1}(A)\subset L^{p_0}(A)\subset L^1(A)$ for $1\leq p_0\leq p_1$. Of course these inclusions are still valid if we consider $A=\mathbb{S}^{n-1}$.

\item For $1\leq p<+\infty$, Lorentz spaces $L^{p,\infty}(A)$ with $A=\mathbb{R}^n$ or $A=\mathbb{S}^{n-1}$ are defined by the condition $\|f\|_ {L^{p,\infty}}=\underset{\lambda>0}{\sup}\{\lambda \times |\{x\in A:|f(x)|> \lambda\}|^{1/p}\}<+\infty$. Recall now that by the real interpolation theory (see \cite[Theorem 5.2.1]{Bergh}) we have for some parameter $0<\theta<1$ the identity
$$(L^p(A), L^\infty(A))_{\theta,\infty}=L^{\frac{p}{1-\theta},\infty}(A).$$
Recall that we always have $L^{\frac{p}{1-\theta},\infty}(A)\subset L^p(A)+L^\infty(A)$, but if the set $A\subset \mathbb{R}^n$ is bounded, we also have the space inclusions
$$L^{\frac{p}{1-\theta},\infty}(A)\subset L^p(A)+L^\infty(A)\subset L^p(A)+L^p(A)=L^p(A).$$ 
Thus, in the particular case of $A=\mathbb{S}^{n-1}$, since $\sigma(\mathbb{S}^{n-1})<+\infty$, we deduce that the Lorentz spaces $L^{q,\infty}(\mathbb{S}^{n-1})$ are embedded in the Lebesgue spaces $L^\rho(\mathbb{S}^{n-1})$ as long as $q>\rho$. In particular, we have 
\begin{equation}\label{InclusionLorentzLpUloc}
L^{n,\infty}(\mathbb{S}^{n-1})\subset L^\rho(\mathbb{S}^{n-1}), \qquad \mbox{if }\quad 1\leq \rho<n.
\end{equation}
\item We consider now the homogeneous Morrey space that are a useful generalization of Lebesgue spaces. Indeed, for $1\leq p\leq q<+\infty$ we define the Morrey space $\dot{\mathcal{M}}^{p,q}(\mathbb{R}^n)$ as the space of measurable functions $f:\mathbb{R}^n\longrightarrow \mathbb{R}$ that are locally in $L^p$ and such that 
\begin{equation}\label{Def_Morrey_space}
\|f\|_{\dot{\mathcal{M}}^{p,q}}= \underset{x\in \mathbb{R}^n, \; r>0}{\sup}\left(\frac{1}{|B(x,r)|^{(1-\frac{p}{q})}}\int_{B(x,r)}|f(y)|^p dy\right)^{\frac{1}{p}}<+\infty.
\end{equation}
Of course, if $p=q$ we have $\dot{\mathcal{M}}^{p,p}(\mathbb{R}^n)=L^p(\mathbb{R}^n)$, however if $1\leq p_0 \leq p_1 \leq q$ and $1<q<+\infty$, then we have the following space inclusions:
$$L^q(\mathbb{R}^n)=\dot{\mathcal{M}}^{q,q}(\mathbb{R}^n) \subset \dot{\mathcal{M}}^{p_1,q}(\mathbb{R}^n) \subset \dot{\mathcal{M}}^{p_0,q}(\mathbb{R}^n).$$
Note now that for $\rho>0$ we have the identity
\begin{equation}\label{StabilityMorrey}
\||f|^\rho\|_{\dot{\mathcal{M}}^{p,q}}=\|f\|_{\dot{\mathcal{M}}^{\rho p, \rho q}}^\rho.
\end{equation}
\end{itemize}
\noindent To end this section we need to recall an important inequality. Indeed, for a function $f\in \mathcal{C}^\infty_0(\mathbb{R}^n)$ and for all ball $B(x,r)$ such that $B(x,r)\subset supp(f)$ we have the following \emph{Poincaré-Sobolev inequality}:
\begin{equation}\label{PoincareSobolev_inequality}
\left(\frac{1}{|B(x,r)|}\int_{B(x,r)}|f(y)-f_{B_r}|^{q}dy\right)^{\frac{1}{q}}\leq Cr\left(\frac{1}{|B(x,r)|}\int_{B(x,r)}|\nabla f(y)|^{\mathfrak{p}}dy\right)^{\frac{1}{\mathfrak{p}}}. 
\end{equation}
for $1\leq \mathfrak{p}<n$ and $1\leq q\leq \frac{n\mathfrak{p}}{n-\mathfrak{p}}$. 
See a proof of this inequality in \cite[Theorem 3.14]{Kinnunen}. 
\mysection{Proof of the Theorem \ref{Theo_Pointwise}}\label{Secc_ProofTh1}
In this section we prove Theorem \ref{Theo_Pointwise}. Let us start by defining 
\begin{equation}\label{Def_OperatorTstar}
T^*_{\Omega}(f)(x)=\underset{t>0}{\sup}\left|\int_{\{|y|>t\}}\frac{\Omega(y/|y|)}{|y|^n}f(x-y)dy\right|,
\end{equation}
and we consider 
$$T^t_{\Omega}(f)(x)=\int_{\{|y|>t\}}\frac{\Omega(y/|y|)}{|y|^n}f(x-y)dy,$$
note that we have $T^*_{\Omega}(f)(x)=\underset{t>0}{\sup}|T^t_{\Omega}(f)(x)|$ and that $|T_\Omega(f)(x)|\leq T^*_{\Omega}(f)(x)$.\\

Now, for a function $f\in \mathcal{C}^\infty_0(\mathbb{R}^n)$ and for some $k_0\in \mathbb{Z}$ so that $2^{k_0-2}<t \leq 2^{k_0-1}$, we write
$$T^t_{\Omega}(f)(x)=\int_{\{t<|y|\leq 2^{k_0-1}\}}\frac{\Omega(y/|y|)}{|y|^n}f(x-y)dy+\sum_{k\geq k_0}\int_{\{2^{k-1}<|y|\leq 2^{k}\}}\frac{\Omega(y/|y|)}{|y|^n}f(x-y)dy.$$
Using the fact that the function $\Omega$ is of null integral, we can introduce some constants in the previous expression to obtain
$$T^t_{\Omega}(f)(x)=\int_{\{t<|y|\leq 2^{k_0-1}\}}\frac{\Omega(y/|y|)}{|y|^n}(f(x-y)-c_{k_0})dy+\sum_{k\geq k_0}\int_{\{2^{k-1}<|y|\leq 2^{k}\}}\frac{\Omega(y/|y|)}{|y|^n}(f(x-y)-c_k)dy,$$
from which we deduce
\begin{eqnarray*}
|T^t_{\Omega}(f)(x)|&\leq &\sum_{k\in \mathbb{Z}}\int_{\{2^{k-1}<|y|\leq 2^{k}\}}\left|\frac{\Omega(y/|y|)}{|y|^n}(f(x-y)-c_k)\right|dy\\
&\leq &C\sum_{k\in \mathbb{Z}}\frac{1}{2^{kn}}\int_{\{|y|\leq 2^{k}\}}\left|\Omega(y/|y|)(f(x-y)-c_k)\right|dy.
\end{eqnarray*}
Now, by the H\"older inequality with $\frac{1}{\rho}+\frac{1}{\rho'}=1$ and $1<\rho<n$, we write
$$|T^t_{\Omega}(f)(x)|\leq \sum_{k\in \mathbb{Z}}\frac{1}{2^{kn}}\left(\int_{\{|y|\leq 2^{k}\}}|\Omega(y/|y|)|^\rho dy\right)^{\frac{1}{\rho}}\left(\int_{\{|y|\leq 2^{k}\}}|f(x-y)-c_k|^{\rho'}dy\right)^{\frac{1}{\rho'}}.$$
Introducing the variable $z=2^{-k} y$, by a change of variables in the first integral above we obtain
$$|T^t_{\Omega}(f)(x)|\leq C\sum_{k\in \mathbb{Z}}\frac{1}{2^{kn(1-\frac{1}{\rho})}}\left(\int_{\{|z|\leq 1\}}|\Omega(z/|z|)|^\rho dz\right)^{\frac{1}{\rho}}\left(\int_{\{|y|\leq 2^{k}\}}|f(x-y)-c_k|^{\rho'}dy\right)^{\frac{1}{\rho'}},$$
and rewriting this formula we have
$$|T^t_{\Omega}(f)(x)|\leq C\sum_{k\in \mathbb{Z}}\left(\int_{\{|z|\leq 1\}}|\Omega(z/|z|)|^\rho dz\right)^{\frac{1}{\rho}}\left(\frac{1}{2^{kn}}\int_{\{|y|\leq 2^{k}\}}|f(x-y)-c_k|^{\rho'}dy\right)^{\frac{1}{\rho'}}.$$
For the second integral above, we consider now the ball $B(x,2^k)$ and we fix the constant $c_k=f_{B_k}=\frac{1}{|B(x,2^k)|}\displaystyle{\int_{B(x,2^k)}f(y)dy}$, so we can write (since $\omega_n 2^{kn}=|B(x,2^k)|$, where $\omega_n=|B(0,1)|$ is the volume of the $n$-dimensional unit ball):
$$|T^t_{\Omega}(f)(x)|\leq C\sum_{k\in \mathbb{Z}}\left(\int_{\{|z|\leq 1\}}|\Omega(z/|z|)|^\rho dz\right)^{\frac{1}{\rho}}\left(\frac{1}{|B(x,2^k)|}\int_{B(x,2^k)}|f(y)-f_{B_k}|^{\rho'}dy\right)^{\frac{1}{\rho'}}.$$
We study now more in detail the first integral above, we thus have
\begin{eqnarray*}
|T^t_{\Omega}(f)(x)|&\leq &C\left(\int_{0}^1\int_{\mathbb{S}^{n-1}}|\Omega(\xi/|\xi|)|^\rho d\sigma(\xi)r^{n-1}dr\right)^{\frac{1}{\rho}}\sum_{k\in \mathbb{Z}}\left(\frac{1}{|B(x,2^k)|}\int_{B(x,2^k)}|f(y)-f_{B_k}|^{\rho'}dy\right)^{\frac{1}{\rho'}}\\
&\leq &C\left(\int_{\mathbb{S}^{n-1}}|\Omega(\xi/|\xi|)|^\rho d\sigma(\xi)\right)^{\frac{1}{\rho}}\sum_{k\in \mathbb{Z}}\left(\frac{1}{|B(x,2^k)|}\int_{B(x,2^k)}|f(y)-f_{B_k}|^{\rho'}dy\right)^{\frac{1}{\rho'}},\end{eqnarray*}
so we obtain
$$|T^t_{\Omega}(f)(x)|\leq C \|\Omega\|_{L^\rho(\mathbb{S}^{n-1})}\sum_{k\in \mathbb{Z}}\left(\frac{1}{|B(x,2^k)|}\int_{B(x,2^k)}|f(y)-f_{B_k}|^{\rho'}dy\right)^{\frac{1}{\rho'}}.$$
\begin{Remarque}
Note that since $1< \rho<n$, by (\ref{InclusionLorentzLpUloc}) we have $L^{n,\infty}(\mathbb{S}^{n-1})\subset L^\rho(\mathbb{S}^{n-1})$ and thus the norm $ \|\Omega\|_{L^\rho(\mathbb{S}^{n-1})}$ induces a refinement with respect to the norm $ \|\Omega\|_{L^{n,\infty}}$ used in \cite{Hoang}. Note also that if $\rho<n$ then we have $\frac{n}{n-1}<\rho'$ as we have $\frac{1}{\rho}+\frac{1}{\rho'}=1$.
\end{Remarque}
Now we apply the Poincaré-Sobolev inequality given in (\ref{PoincareSobolev_inequality}) to obtain
\begin{eqnarray}
|T^t_{\Omega}(f)(x)|&\leq &C \|\Omega\|_{L^\rho(\mathbb{S}^{n-1})}\sum_{k\in \mathbb{Z}}\left(\frac{1}{|B(x,2^k)|}\int_{B(x,2^k)}|f(y)-f_{B_k}|^{\rho'}dy\right)^{\frac{1}{\rho'}}\notag\\
&\leq &C \|\Omega\|_{L^\rho(\mathbb{S}^{n-1})}\sum_{k\in \mathbb{Z}}
2^k\left(\frac{1}{|B(x,2^k)|}\int_{B(x,2^k)}|\nabla f(y)|^{\mathfrak{p}}dy\right)^{\frac{1}{\mathfrak{p}}},\label{Somme1}
\end{eqnarray}
where $\frac{n}{n-1}<\rho'$ (since $1<\rho<n$ and $\frac{1}{\rho}+\frac{1}{\rho'}=1$) and where $\rho'\leq \frac{n\mathfrak{p}}{n-\mathfrak{p}}$. Note that we thus have $\frac{n}{n-1}<\rho'\leq \frac{n\mathfrak{p}}{n-\mathfrak{p}}$ which leads us to the condition $1<\frac{\rho n}{\rho n+\rho-n}=\overline{\rho}\leq \mathfrak{p}<n$. 
We study now the sum
$$\mathcal{S}=\sum_{k\in \mathbb{Z}}2^k\left(\frac{1}{|B(x,2^k)|}\int_{B(x,2^k)}|\nabla f(y)|^{\mathfrak{p}}dy\right)^{\frac{1}{\mathfrak{p}}},$$
and since we can use the following triangle inequality, for all fixed $x\in \mathbb{R}^n$ and $k\in \mathbb{Z}$:
\begin{eqnarray*}
\left(\frac{1}{|B(x,2^k)|}\int_{B(x,2^k)}|\nabla f(y)|^{\mathfrak{p}}dy\right)^{\frac{1}{\mathfrak{p}}}&\leq &\left(\frac{1}{|B(x,2^k)|}\int_{B(x,2^k)}|\nabla f(y)|^{\mathfrak{p}}\mathds{1}_{\{|x-y|\leq 2^{k-1}\}}dy\right)^{\frac{1}{\mathfrak{p}}}\\
&&+\left(\frac{1}{|B(x,2^k)|}\int_{B(x,2^k)}|\nabla f(y)|^{\mathfrak{p}}\mathds{1}_{\{2^{k-1}\leq |x-y|\leq 2^k\}}dy\right)^{\frac{1}{\mathfrak{p}}},
\end{eqnarray*}
we then have
\begin{equation}\label{Somme2}
\begin{split}
\mathcal{S}&\leq \underbrace{\sum_{k\in \mathbb{Z}}2^k\left(\frac{1}{|B(x,2^k)|}\int_{\{|x-y|\leq 2^{k-1}\}}|\nabla f(y)|^{\mathfrak{p}}dy\right)^{\frac{1}{\mathfrak{p}}}}_{\mathcal{S}_1} \\
&+\underbrace{\sum_{k\in \mathbb{Z}}2^k\left(\frac{1}{|B(x,2^k)|}\int_{\{2^{k-1}\leq |x-y|\leq 2^k\}}|\nabla f(y)|^{\mathfrak{p}}dy\right)^{\frac{1}{\mathfrak{p}}}}_{\mathcal{S}_2}.
\end{split}
\end{equation}
\begin{itemize}
\item For the first term of (\ref{Somme2}) we have
\begin{eqnarray*}
\mathcal{S}_1&=&\sum_{k\in \mathbb{Z}}2^k\left(\frac{1}{|B(x,2^k)|}\int_{\{|x-y|\leq 2^{k-1}\}}|\nabla f(y)|^{\mathfrak{p}}dy\right)^{\frac{1}{\mathfrak{p}}}\\
&=&2^{1-\frac{n}{\mathfrak{p}}}\sum_{k\in \mathbb{Z}}2^{k-1}\left(\frac{1}{|B(x,2^{k-1})|}\int_{\{|x-y|\leq 2^{k-1}\}}|\nabla f(y)|^{\mathfrak{p}}dy\right)^{\frac{1}{\mathfrak{p}}},
\end{eqnarray*}
which is 
$$\mathcal{S}_1=2^{1-\frac{n}{\mathfrak{p}}}\sum_{k\in \mathbb{Z}}2^k\left(\frac{1}{|B(x,2^k)|}\int_{B(x,2^k)}|\nabla f(y)|^{\mathfrak{p}}dy\right)^{\frac{1}{\mathfrak{p}}},$$
and we obtain the formula
\begin{equation}\label{FormulaS1}
\mathcal{S}_1=2^{1-\frac{n}{\mathfrak{p}}}\mathcal{S}.
\end{equation}
\item For the second term of (\ref{Somme2}) we have
$$\mathcal{S}_2=\sum_{k\in \mathbb{Z}}2^k\left(\frac{1}{|B(x,2^k)|}\int_{\{2^{k-1}\leq |x-y|\leq 2^k\}}|\nabla f(y)|^{\mathfrak{p}}dy\right)^{\frac{1}{\mathfrak{p}}}.$$
Introducing a parameter $0<K<+\infty$ that will be fixed later, we write
\begin{eqnarray}
\mathcal{S}_2&=&\underbrace{\sum_{k\leq  \lfloor \log_2(K)\rfloor }2^k\left(\frac{1}{|B(x,2^k)|}\int_{\{2^{k-1}\leq |x-y|\leq 2^k\}}|\nabla f(y)|^{\mathfrak{p}}dy\right)^{\frac{1}{\mathfrak{p}}}}_{(A)}\notag\\
&&+\underbrace{\sum_{k>  \lfloor \log_2(K)\rfloor }2^k\left(\frac{1}{|B(x,2^k)|}\int_{\{2^{k-1}\leq |x-y|\leq 2^k\}}|\nabla f(y)|^{\mathfrak{p}}dy\right)^{\frac{1}{\mathfrak{p}}}}_{(B)}.\label{Somme22}
\end{eqnarray}
The term $(A)$ in the expression (\ref{Somme22}) above is treated as follows: 
\begin{eqnarray*}
(A)&=&\sum_{k\leq  \lfloor \log_2(K)\rfloor }2^k\left(\frac{1}{|B(x,2^k)|}\int_{\{2^{k-1}\leq |x-y|\leq 2^k\}}|\nabla f(y)|^{\mathfrak{p}}dy\right)^{\frac{1}{\mathfrak{p}}}\\
&\leq &\sum_{k\leq  \lfloor \log_2(K)\rfloor }2^k\left(\frac{1}{|B(x,2^k)|}\int_{B(x,2^k)}|\nabla f(y)|^{\mathfrak{p}}dy\right)^{\frac{1}{\mathfrak{p}}},
\end{eqnarray*}
and since we have the control $\displaystyle{\left(\frac{1}{|B(x,2^k)|}\int_{B(x,2^k)} |\nabla f(y)|^{\mathfrak{p}} dy\right)^{\frac{1}{\mathfrak{p}}}\leq \big(\mathscr{M}(|\nabla f|^{\mathfrak{p}})(x)\big)^{\frac{1}{\mathfrak{p}}}}$, we obtain
\begin{eqnarray}
(A)&\leq &\sum_{k\leq  \lfloor \log_2(K)\rfloor }2^k\big(\mathscr{M}(|\nabla f|^{\mathfrak{p}})(x)\big)^{\frac{1}{\mathfrak{p}}}=\big(\mathscr{M}(|\nabla f|^{\mathfrak{p}})(x)\big)^{\frac{1}{\mathfrak{p}}} \sum_{k\leq  \lfloor \log_2(K)\rfloor }2^k\notag \\
&\leq& C K\big(\mathscr{M}(|\nabla f|^{\mathfrak{p}})(x)\big)^{\frac{1}{\mathfrak{p}}}.\label{SommeA}
\end{eqnarray}
For the term $(B)$ in the expression (\ref{Somme22}) we write
\begin{eqnarray*}
(B)&=&\sum_{k>  \lfloor \log_2(K)\rfloor }2^k\left(\frac{1}{|B(x,2^k)|}\int_{\{2^{k-1}\leq |x-y|\leq 2^k\})}|\nabla f(y)|^{\mathfrak{p}}dy\right)^{\frac{1}{\mathfrak{p}}}\\
&\leq &\sum_{k>  \lfloor \log_2(K)\rfloor }2^k\left(\frac{1}{|B(x,2^k)|}\int_{B(x,2^k)}|\nabla f(y)|^{\mathfrak{p}}dy\right)^{\frac{1}{\mathfrak{p}}},
\end{eqnarray*}
and we have 
\begin{eqnarray*}
(B)&\leq &\sum_{k>  \lfloor \log_2(K)\rfloor }2^k\left(\frac{|B(x,2^k)|^{-\frac{\mathfrak{p}}{\mathfrak{q}}}}{|B(x,2^k)|^{1-\frac{\mathfrak{p}}{\mathfrak{q}}}}\int_{B(x,2^k)}|\nabla f(y)|^{\mathfrak{p}}dy\right)^{\frac{1}{\mathfrak{p}}}\\
&\leq &C\sum_{k>  \lfloor \log_2(K)\rfloor }2^{k(1-\frac{n}{\mathfrak{q}})}\left(\frac{1}{|B(x,2^k)|^{(1-\frac{\mathfrak{p}}{\mathfrak{q}})}}\int_{B(x,2^k)}|\nabla f(y)|^{\mathfrak{p}}dy\right)^{\frac{1}{\mathfrak{p}}}.
\end{eqnarray*}
Using the definition of the Morrey spaces $\dot{\mathcal{M}}^{\mathfrak{p},\mathfrak{q}}(\mathbb{R}^n)$ given in the formula (\ref{Def_Morrey_space}), we obtain
\begin{eqnarray*}
(B)\leq C\|\nabla f\|_{\dot{\mathcal{M}}^{\mathfrak{p},\mathfrak{q}}}\sum_{k>  \lfloor \log_2(K)\rfloor }2^{k(1-\frac{n}{\mathfrak{q}})}.
\end{eqnarray*}
But since $1-\frac{n}{\mathfrak{q}}<0$, as we have by hypothesis $\mathfrak{q}<n$ (recall that we have $1<\frac{\rho n}{\rho n+\rho-n}=\overline{\rho}\leq \mathfrak{p}\leq \mathfrak{q}<n$), we obtain that the previous sum converges and we can write 
\begin{equation}\label{SommeB}
(B)\leq C\|\nabla f\|_{\dot{\mathcal{M}}^{\mathfrak{p},\mathfrak{q}}}K^{1-\frac{n}{\mathfrak{q}}}.
\end{equation}
With the estimates (\ref{SommeA}) and (\ref{SommeB}) at hand, we come back to the expression (\ref{Somme22}) to obtain the inequality
$$\mathcal{S}_2\leq C\left(K\big(\mathscr{M}(|\nabla f|^{\mathfrak{p}})(x)\big)^{\frac{1}{\mathfrak{p}}}+\|\nabla f\|_{\dot{\mathcal{M}}^{\mathfrak{p},\mathfrak{q}}}K^{1-\frac{n}{\mathfrak{q}}}\right).$$
If we set $K=\left(\frac{\|\nabla f\|_{\dot{\mathcal{M}}^{\mathfrak{p},\mathfrak{q}}}}{\big(\mathscr{M}(|\nabla f|^{\mathfrak{p}})(x)\big)^{\frac{1}{\mathfrak{p}}}}\right)^{\frac{\mathfrak{q}}{n}}$, we have
\begin{equation}\label{FormulaS2} 
\mathcal{S}_2\leq C \big(\mathscr{M}(|\nabla f|^{\mathfrak{p}})(x)\big)^{\frac{1}{\mathfrak{p}}-\frac{\mathfrak{q}}{\mathfrak{p}n}}\|\nabla f\|_{\dot{\mathcal{M}}^{\mathfrak{p},\mathfrak{q}}}^{\frac{\mathfrak{q}}{n}}.
\end{equation}
\end{itemize}
With the estimate (\ref{FormulaS1}) for the term $\mathcal{S}_1$ and the inequality (\ref{FormulaS2}) for the term $\mathcal{S}_2$, getting back to the expression (\ref{Somme2}) we have the control:
\begin{eqnarray*}
\mathcal{S}&\leq & \mathcal{S}_1+\mathcal{S}_2\\
&\leq & 2^{1-\frac{n}{\mathfrak{p}}}\mathcal{S}+C \big(\mathscr{M}(|\nabla f|^{\mathfrak{p}})(x)\big)^{\frac{1}{\mathfrak{p}}-\frac{\mathfrak{q}}{\mathfrak{p}n}}\|\nabla f\|_{\dot{\mathcal{M}}^{\mathfrak{p},\mathfrak{q}}}^{\frac{\mathfrak{q}}{n}}.
\end{eqnarray*}
Since $\mathfrak{p}<n$ we have $2^{1-\frac{n}{\mathfrak{p}}}<1$ and we obtain 
$$\mathcal{S}(1- 2^{1-\frac{n}{\mathfrak{p}}})\leq C \big(\mathscr{M}(|\nabla f|^{\mathfrak{p}})(x)\big)^{\frac{1}{\mathfrak{p}}-\frac{\mathfrak{q}}{\mathfrak{p}n}}\|\nabla f\|_{\dot{\mathcal{M}}^{\mathfrak{p},\mathfrak{q}}}^{\frac{\mathfrak{q}}{n}},$$
from which we deduce that $\mathcal{S}\leq C\big(\mathscr{M}(|\nabla f|^{\mathfrak{p}})(x)\big)^{\frac{1}{\mathfrak{p}}-\frac{\mathfrak{q}}{\mathfrak{p}n}}\|\nabla f\|_{\dot{\mathcal{M}}^{\mathfrak{p},\mathfrak{q}}}^{\frac{\mathfrak{q}}{n}}$. Thus, coming back to the expression (\ref{Somme1}), we have 
$$|T^t_{\Omega}(f)(x)|\leq C \|\Omega\|_{L^\rho(\mathbb{S}^{n-1})}\big(\mathscr{M}(|\nabla f|^{\mathfrak{p}})(x)\big)^{\frac{1}{\mathfrak{p}}-\frac{\mathfrak{q}}{\mathfrak{p}n}}\|\nabla f\|_{\dot{\mathcal{M}}^{\mathfrak{p},\mathfrak{q}}}^{\frac{\mathfrak{q}}{n}},$$
and from this uniform estimate we can obtain the control
\begin{equation}\label{EstimationInitialRefined1}
T^*_{\Omega}(f)(x)\leq C \|\Omega\|_{L^\rho(\mathbb{S}^{n-1})}\big(\mathscr{M}(|\nabla f|^{\mathfrak{p}})(x)\big)^{\frac{1}{\mathfrak{p}}-\frac{\mathfrak{q}}{\mathfrak{p}n}}\|\nabla f\|_{\dot{\mathcal{M}}^{\mathfrak{p},\mathfrak{q}}}^{\frac{\mathfrak{q}}{n}},
\end{equation}
and this finishes the proof of the Theorem \ref{Theo_Pointwise} since we have $|T_\Omega(f)(x)|\leq T^*_{\Omega}(f)(x)$. \hfill $\blacksquare$
\mysection{Proof of the Theorem \ref{Theo_RefinedSobolevInequalities}}\label{Secc_ProofTheo2}
Our starting point is the pointwise estimate obtained in the previous Theorem \ref{Theo_Pointwise}:
\begin{equation}\label{Modified_hedbergEstimate}
|T_\Omega(f)(x)|\leq C_\Omega \big(\mathscr{M}(|\nabla f|^{\mathfrak{p}})(x)\big)^{\frac{1}{\mathfrak{p}}-\frac{\mathfrak{q}}{\mathfrak{p}n}}\|\nabla f\|_{\dot{\mathcal{M}}^{\mathfrak{p},\mathfrak{q}}}^{\frac{\mathfrak{q}}{n}},
\end{equation}
Thus, taking the $L^r$-norm to both sides of the estimate (\ref{Modified_hedbergEstimate}) we have the inequality
$$\|T_\Omega(f)\|_{L^r}\leq C_\Omega  \left\|\big(\mathscr{M}(|\nabla f|^{\mathfrak{p}})\big)^{\frac{1}{\mathfrak{p}}-\frac{\mathfrak{q}}{\mathfrak{p}n}}\right\|_{L^r}\|\nabla f\|_{\dot{\mathcal{M}}^{\mathfrak{p},\mathfrak{q}}}^{\frac{\mathfrak{q}}{n}},$$
which we rewrite as (using the property $\||f|^s\|_{L^p}=\|f\|_{L^{s p}}^{s}$ for the Lebesgue norms):
$$\|T_\Omega(f)\|_{L^r}\leq C_\Omega \left\|\mathscr{M}(|\nabla f|^{\mathfrak{p}})\right\|_{L^{r(\frac{1}{\mathfrak{p}}-\frac{\mathfrak{q}}{\mathfrak{p}n})}}^{\frac{1}{\mathfrak{p}}-\frac{\mathfrak{q}}{\mathfrak{p}n}}\|\nabla f\|_{\dot{\mathcal{M}}^{\mathfrak{p},\mathfrak{q}}}^{\frac{\mathfrak{q}}{n}}.$$
Since by hypothesis we have $r=\frac{n\sigma}{n-\mathfrak{q}}>1$ and $\sigma>\mathfrak{p}$, 
we deduce that $r(\frac{1}{\mathfrak{p}}-\frac{\mathfrak{q}}{\mathfrak{p}n})=\frac{r}{\mathfrak{p}}(\frac{n-\mathfrak{q}}{n})=\frac{\sigma}{\mathfrak{p}}>1$, thus the maximal function $\mathscr{M}$ is bounded in the Lebesgue space $L^{r(\frac{1}{\mathfrak{p}}-\frac{\mathfrak{q}}{\mathfrak{p}n})}(\mathbb{R}^n)$ and we can write
$$\|T_\Omega(f)\|_{L^r}\leq C_\Omega \left\||\nabla f|^{\mathfrak{p}}\right\|_{L^{r(\frac{1}{\mathfrak{p}}-\frac{\mathfrak{q}}{\mathfrak{p}n})}}^{\frac{1}{\mathfrak{p}}-\frac{\mathfrak{q}}{\mathfrak{p}n}}\|\nabla f\|_{\dot{\mathcal{M}}^{\mathfrak{p},\mathfrak{q}}}^{\frac{\mathfrak{q}}{n}}\leq C_\Omega \left\|\nabla f\right\|_{L^{r(1-\frac{\mathfrak{q}}{n})}}^{1-\frac{\mathfrak{q}}{n}}\|\nabla f\|_{\dot{\mathcal{M}}^{\mathfrak{p},\mathfrak{q}}}^{\frac{\mathfrak{q}}{n}},$$
where we used again the property $\||f|^s\|_{L^p}=\|f\|_{L^{s p}}^{s}$. Since $\theta=\frac{\mathfrak{q}}{n}<1$ and since $\sigma=r(1-\frac{\mathfrak{q}}{n})$, we can finally obtain the inequality 
$$\|T_\Omega(f)\|_{L^r}\leq C_\Omega \left\|\nabla f\right\|_{L^{\sigma}}^{1-\theta}\|\nabla f\|_{\dot{\mathcal{M}}^{\mathfrak{p},\mathfrak{q}}}^{\theta},$$
and the Theorem \ref{Theo_RefinedSobolevInequalities} is now proven. \hfill$\blacksquare$
\begin{Remarque}
Note that if $\mathfrak{p}<\sigma=\mathfrak{q}<n$, since we have $L^{\mathfrak{q}}(\mathbb{R}^n)=\dot{\mathcal{M}}^{\mathfrak{q},\mathfrak{q}}(\mathbb{R}^n)\subset \dot{\mathcal{M}}^{\mathfrak{p},\mathfrak{q}}(\mathbb{R}^n)$, we obtain
$$\|T_\Omega(f)\|_{L^r}\leq C_\Omega \left\|\nabla f\right\|_{L^{\mathfrak{q}}}^{1-\theta}\|\nabla f\|_{\dot{\mathcal{M}}^{\mathfrak{p},\mathfrak{q}}}^{\theta}\leq C_\Omega \left\|\nabla f\right\|_{L^{\mathfrak{q}}}^{1-\theta}\|\nabla f\|_{\dot{\mathcal{M}}^{\mathfrak{q},\mathfrak{q}}}^{\theta}=C_\Omega \left\|\nabla f\right\|_{L^{\mathfrak{q}}},$$
which is a new version of the ``classical'' Sobolev inequalities for rough type operators $T_\Omega$ associated to a function $\Omega \in L^\rho(\mathbb{S}^{n-1})$ where $1<\rho<n$, $\frac{\rho n}{\rho n+\rho-n} <\mathfrak{q}<n$ and $r=\frac{n\mathfrak{q}}{n-\mathfrak{q}}>1$.
\end{Remarque}
\mysection{Weighted inequalities}\label{Secc_Weights}
The inequality (\ref{InequalitySobolev_Intro}) can be easily generalized by considering suitable weights. First, recall that for a generic weight $w:\mathbb{R}^n\longrightarrow \mathbb{R}^+$, for $1\leq p<+\infty$ we define the weighted Lebesgue spaces $L^p(w)$ by the condition
\begin{equation}\label{Def_WeightedLebesgue}
\|f\|_{L^p(w)}=\left(\int_{\mathbb{R}^n}|f(x)|^pw(x)dx\right)^{\frac{1}{p}}<+\infty.
\end{equation}
Note that from this definition we deduce, for some $s>0$ the property
\begin{equation}\label{Puissance_WeightedLebesgue}
\||f|^s\|_{L^p(w)}=\|f\|_{L^{s p}(w)}^s.
\end{equation}
Although many type of weights are available in the literature, as we will need to deal at some point with the Hardy-Littlewood maximal function $\mathscr{M}$, it is quite natural to consider weights in the $A_\delta$ class: for $1<\delta<+\infty$ we will say that a weight $w$ belongs to the $A_\delta$ class if $w^{-1}$ is locally integrable and if 
$$[w]_{\delta}=\underset{B}{\sup}\left(\frac{1}{|B|}\int_{B}w(x)dx\right)\left(\frac{1}{|B|}\int_{B}w(x)^{-\frac{1}{\delta-1}}dx\right)^{\delta-1}<+\infty.$$
Note that the $A_\delta$ class gives a quite natural framework to obtain the following estimate
\begin{equation}\label{BoundednessMaximale}
\|\mathscr{M}(f)\|_{L^\delta(w)}\leq C\|f\|_{L^\delta(w)},
\end{equation}
and this boundedness property is actually equivalent to the fact that $w\in A_\delta$. See the book \cite{Grafakos} for more details and properties of this class of weights.\\

In this context we have the following result:
\begin{Corollaire}[\bf One weighted inequality]\label{Coro_Weighted1}
Over the space $\mathbb{R}^n$ with $n\geq 2$, consider $\Omega$ a function such that $\Omega \in L^1(\mathbb{S}^{n-1})$, $\displaystyle{\int_{\mathbb{S}^{n-1}}\Omega \ d\sigma=0}$ and such that $\Omega\in L^\rho(\mathbb{S}^{n-1})$ with $1<\rho<n$ and consider the operator $T_{\Omega}$ associated to the function $\Omega$ as defined in (\ref{Def_Operator}). \\

\noindent Assume that $f:\mathbb{R}^n\longrightarrow \mathbb{R}$ is a regular function such that $\nabla f\in \dot{\mathcal{M}}^{\mathfrak{p},\mathfrak{q}}(\mathbb{R}^n)$ where $1< \mathfrak{p}\leq \mathfrak{q}<n$ and where $\mathfrak{p}$ is a real parameter such that $1<\overline{\rho}\leq  \mathfrak{p}$ with $\overline{\rho}=\frac{\rho n}{\rho n+\rho-n}$. If we have $\nabla f\in L^{\sigma}(\omega)$ where $\sigma > \mathfrak{p}$ and where $\omega$ is a weight such that $\omega\in A_{\frac{r}{\mathfrak{p}}(1-\theta)}$ where $r=\frac{n\sigma}{n-\mathfrak{q}}>1$ and $\theta=\frac{\mathfrak{q}}{n}$, then we have the inequality
$$\|T_\Omega(f)\|_{L^r(\omega)}\leq C_\Omega \left\|\nabla f\right\|_{L^{\sigma}(\omega)}^{1-\theta}\|\nabla f\|_{\dot{\mathcal{M}}^{\mathfrak{p},\mathfrak{q}}}^{\theta}.$$
\end{Corollaire}
\noindent {\bf Proof.} By the pointwise estimate (\ref{Modified_hedbergEstimate}) we have
$$|T_\Omega(f)(x)|\leq C_\Omega \big(\mathscr{M}(|\nabla f|^{\mathfrak{p}})(x)\big)^{\frac{1}{\mathfrak{p}}-\frac{\mathfrak{q}}{\mathfrak{p}n}}\|\nabla f\|_{\dot{\mathcal{M}}^{\mathfrak{p},\mathfrak{q}}}^{\frac{\mathfrak{q}}{n}},$$
which we raise to the power $r$, multiply by the weight $w$ and integrate to obtain
$$\|T_\Omega(f)\|_{L^r(w)}\leq C_\Omega \left\|\big(\mathscr{M}(|\nabla f|^{\mathfrak{p}})\big)^{\frac{1}{\mathfrak{p}}-\frac{\mathfrak{q}}{\mathfrak{p}n}}\right\|_{L^r(w)}\|\nabla f\|_{\dot{\mathcal{M}}^{\mathfrak{p},\mathfrak{q}}}^{\frac{\mathfrak{q}}{n}}.$$
By the property (\ref{Puissance_WeightedLebesgue}) of the weighted Lebesgue spaces we have 
$$\|T_\Omega(f)\|_{L^r(w)}\leq C_\Omega  \left\|\mathscr{M}(|\nabla f|^{\mathfrak{p}})\right\|_{L^{r(\frac{1}{\mathfrak{p}}-\frac{\mathfrak{q}}{\mathfrak{p}n})}(w)}^{\frac{1}{\mathfrak{p}}-\frac{\mathfrak{q}}{\mathfrak{p}n}}\|\nabla f\|_{\dot{\mathcal{M}}^{\mathfrak{p},\mathfrak{q}}}^{\frac{\mathfrak{q}}{n}}.$$
Since the weight $\omega$ belongs to the class $A_{\frac{r}{\mathfrak{p}}(1-\theta)}$ (recall that we have $\frac{r}{\mathfrak{p}}(1-\theta)=r(\frac{1}{\mathfrak{p}}-\frac{\mathfrak{q}}{\mathfrak{p}n})>1$) then the maximal function $\mathscr{M}$ is bounded in the Lebesgue spaces $L^{r(\frac{1}{\mathfrak{p}}-\frac{\mathfrak{q}}{\mathfrak{p}n})}(w)$, and we obtain
\begin{eqnarray*}
\|T_\Omega(f)\|_{L^r(w)}&\leq& C_\Omega \left\||\nabla f|^{\mathfrak{p}}\right\|_{L^{r(\frac{1}{\mathfrak{p}}-\frac{\mathfrak{q}}{\mathfrak{p}n})}(w)}^{\frac{1}{\mathfrak{p}}-\frac{\mathfrak{q}}{\mathfrak{p}n}}\|\nabla f\|_{\dot{\mathcal{M}}^{\mathfrak{p},\mathfrak{q}}}^{\frac{\mathfrak{q}}{n}}\\[2mm]
&\leq & C_\Omega \|\nabla f\|_{L^{r(1-\frac{\mathfrak{q}}{n})}(w)}^{1-\frac{\mathfrak{q}}{n}}\|\nabla f\|_{\dot{\mathcal{M}}^{\mathfrak{p},\mathfrak{q}}}^{\frac{\mathfrak{q}}{n}},
\end{eqnarray*}
where we applied again the property (\ref{Puissance_WeightedLebesgue}). Recalling that $\theta=\frac{\mathfrak{q}}{n}$, we have $\sigma=r(1-\theta)$ and we can write
$$\|T_\Omega(f)\|_{L^r(\omega)}\leq C_\Omega \|\nabla f\|_{L^{\sigma}(\omega)}^{1-\theta}\|\nabla f\|_{\dot{\mathcal{M}}^{\mathfrak{p},\mathfrak{q}}}^{\theta},$$
which is the announced estimate. The proof of the corollary is now complete.  \hfill $\blacksquare$
\mysection{Refined Sobolev inequalities in the framework of Orlicz spaces}\label{Secc_Orlicz}
We will consider here an extension of the Theorem \ref{Theo_RefinedSobolevInequalities} to the framework of \emph{Orlicz spaces}. To present these spaces, we first recall that if $a:[0,+\infty[\longrightarrow[0,+\infty[$ is a left-continuous non decreasing function with $a(0)=0$, we can consider the corresponding \emph{Young function} $A(t)=\displaystyle{\int_{0}^ta(s)ds}$. The Orlicz space $L^A(\mathbb{R}^n)$ associated to the Young function $A$ is then defined as the set of measurable functions $f:\mathbb{R}^n\longrightarrow \mathbb{R}$ such that the following Luxemburg norm
\begin{equation}\label{Def_NormOrlicz}
\|f\|_{L^A(\mathbb{R}^n)}=\inf\left\{\lambda > 0: \, \int_{\mathbb{R}^n}A(|f(x)|/\lambda)dx\leq1\right\},
\end{equation}
is finite. Of course we can easily see here that if $A(t)=t^p$ for $1\leq p<+\infty$, we recover the classical Lebesgue spaces. Since the quantity $\|\cdot\|_{L^A}$ is a norm, we have some nice properties: for example, if $f,g$ are two measurable functions such that $|f|\leq |g|$ a.e., then we have the order-reserving property
$$\|f\|_{L^A}\leq \|g\|_{L^A}.$$
However, the Orlicz spaces given by (\ref{Def_NormOrlicz}) with a generic Young function $A$ are too general for our purposes as we need some structure to perform our computations. First we will need the following \emph{rescaling property} as defined in Section 3 of \cite{RaSam}: for any real $s>0$, we define the space $L^A_s(\mathbb{R}^n)$ by the condition
$$L^A_s(\mathbb{R}^n)=\{f:\mathbb{R}^n\longrightarrow \mathbb{R}: \|f\|_{L^A_s(\mathbb{R}^n)}<+\infty\},$$
where 
\begin{equation}\label{Rescal}
\|f\|_{L^A_s}=\inf\left\{\lambda > 0: \, \int_{\mathbb{R}^n}A_s(|f(x)|/\lambda)dx\leq1\right\},
\end{equation}
with $A_s(t)=A(t^s)$. With this definition of the functional $\|\cdot\|_{L^A_s}$ we have the following identity 
\begin{equation}\label{Rescal1}
\||f|^s\|_{L^A}=\|f\|_{L^A_s}^s,
\end{equation}
which will be essential in the sequel. See Lemma 3.2 of \cite{RaSam} for a proof of this fact.\\

Next, it is classical to impose the $\nabla_2$-condition over the Young functions: indeed, a Young function $A$ is said to satisfy the $\nabla_2$-condition, denoted also by $A\in \nabla_2$, if
$$A(r)\leq\frac{1}{2C} A(Cr),\qquad r\geq 0,$$
for some $C > 1$. This condition ensure the boundedness of the Hardy-Littlewood maximal function in the setting of Orlicz spaces: if $A\in  \nabla_2$ we thus have
\begin{equation*}
\|\mathscr{M}(f)\|_{L^{A}}\leq C\|f\|_{L^{A}},
\end{equation*}
see \cite{Cianchi0} for a proof of this fact, see also \cite[Theorem 2]{Derigoz} and the reference there in for more details on the boundedness of the maximal functions in this setting.\\ 

Note that in \cite{Cianchi1} some Sobolev inequalities have been studied in the context of Orlicz spaces. However, and to the best of our knowledge, Sobolev-type inequalities with rough operators seems to be new in this framework. We can thus consider the following result, which is an extension of the Theorem \ref{Theo_RefinedSobolevInequalities} above to the setting of Orlicz spaces:
\begin{Theoreme}
Over the space $\mathbb{R}^n$ with $n\geq 2$, consider $\Omega$ a function such that $\Omega \in L^1(\mathbb{S}^{n-1})$, $\displaystyle{\int_{\mathbb{S}^{n-1}}\Omega \ d\sigma=0}$ and such that $\Omega\in L^\rho(\mathbb{S}^{n-1})$ with $1<\rho<n$ and consider the operator $T_{\Omega}$ associated to the function $\Omega$ as defined in (\ref{Def_Operator}). \\

\noindent Assume that $f:\mathbb{R}^n\longrightarrow \mathbb{R}$ is a regular function such that $\nabla f\in \dot{\mathcal{M}}^{\mathfrak{p},\mathfrak{q}}(\mathbb{R}^n)$ where $1< \mathfrak{p}\leq \mathfrak{q}<n$ and where $\mathfrak{p}$ is a real parameter such that $1<\overline{\rho}\leq \mathfrak{p}$, with $\overline{\rho}=\frac{\rho n}{\rho n+\rho-n}$. \\

\noindent Define $\theta=\frac{\mathfrak{q}}{n}$, if we have $\nabla f\in L^{A}_{(1-\theta)}(\mathbb{R}^n)$ where the space $L^{A}_{(1-\theta)}(\mathbb{R}^n)$ is given by the condition (\ref{Rescal}) and if the Young function $A_{\frac{1}{\mathfrak{p}}(1-\theta)}$ satisfy the $A_2$-condition, then we have the inequality
$$\|T_\Omega(f)\|_{L^A}\leq C_\Omega\|\nabla f\|_{L^{A}_{(1-\theta)}}^{1-\theta}\|\nabla f\|_{\dot{\mathcal{M}}^{\mathfrak{p},\mathfrak{q}}}^{\theta}.$$
\end{Theoreme}
\noindent {\bf Proof.} Once we have at our disposal a good pointwise estimate, the proof is relatively straightforward. Indeed, from the control
$$|T_\Omega(f)(x)|\leq C_\Omega \big(\mathscr{M}(|\nabla f|^{\mathfrak{p}})(x)\big)^{\frac{1}{\mathfrak{p}}-\frac{\mathfrak{q}}{\mathfrak{p}n}}\|\nabla f\|_{\dot{\mathcal{M}}^{\mathfrak{p},\mathfrak{q}}}^{\frac{\mathfrak{q}}{n}},$$
by the order-preserving property of the functional $\|\cdot\|_{L^A}$, we have 
$$\|T_\Omega(f)\|_{L^A}\leq C_\Omega \left\|\mathscr{M}(|\nabla f|^{\mathfrak{p}})^{\frac{1}{\mathfrak{p}}-\frac{\mathfrak{q}}{\mathfrak{p}n}}\right\|_{L^A}\|\nabla f\|_{\dot{\mathcal{M}}^{\mathfrak{p},\mathfrak{q}}}^{\frac{\mathfrak{q}}{n}},$$
now, by the rescaling property (\ref{Rescal1}) we obtain
$$\|T_\Omega(f)\|_{L^A}\leq C_\Omega\left\|\mathscr{M} \left(|\nabla f|^\mathfrak{p}\right)\right\|_{L^{A}_{(\frac{1}{\mathfrak{p}}-\frac{\mathfrak{q}}{\mathfrak{p}n})}}^{(\frac{1}{\mathfrak{p}}-\frac{\mathfrak{q}}{\mathfrak{p}n})}\|\nabla f\|_{\dot{\mathcal{M}}^{\mathfrak{p},\mathfrak{q}}}^{\frac{\mathfrak{q}}{n}}.$$
Since the Young functions $A_{\frac{1}{\mathfrak{p}}(1-\theta)}=A_{(\frac{1}{\mathfrak{p}}-\frac{\mathfrak{q}}{\mathfrak{p}n})}$ satisfy the $\nabla_2$-condition (recall that $\theta=\frac{\mathfrak{q}}{n}$), the Hardy-Littlewood maximal function is bounded in the Orlicz space $L^A_{(\frac{1}{\mathfrak{p}}-\frac{\mathfrak{q}}{\mathfrak{p}n})}$ and thus we can write
$$\|T_\Omega(f)\|_{L^A}\leq C_\Omega\left\||\nabla f|^\mathfrak{p}\right\|_{L^{A}_{(\frac{1}{\mathfrak{p}}-\frac{\mathfrak{q}}{\mathfrak{p}n})}}^{(\frac{1}{\mathfrak{p}}-\frac{\mathfrak{q}}{\mathfrak{p}n})}\|\nabla f\|_{\dot{\mathcal{M}}^{\mathfrak{p},\mathfrak{q}}}^{\frac{\mathfrak{q}}{n}},$$
using again the rescaling property (\ref{Rescal1}) we have
$$\|T_\Omega(f)\|_{L^A}\leq C_\Omega\|\nabla f\|_{L^{A}_{(1-\frac{\mathfrak{q}}{n})}}^{(1-\frac{\mathfrak{q}}{n})}\|\nabla f\|_{\dot{\mathcal{M}}^{\mathfrak{p},\mathfrak{q}}}^{\frac{\mathfrak{q}}{n}},$$
which is the announced inequality since $\theta=\frac{\mathfrak{q}}{n}$. The proof of the theorem is  complete.  \hfill $\blacksquare$
\mysection{Inequalities in classical Lorentz spaces}\label{Secc_ClassicalLorentz}
For $1\leq p<+\infty$ and for $w:\mathbb{R}^+\longrightarrow\mathbb{R}^+$ a weight, we consider there the classical Lorentz space of functions introduced in \cite{Lorentz} and \cite{Lorentz1} defined as
$$\Lambda^p(w)=\left\{f:\|f\|_{\Lambda^p(w)}=\left(\int_0^{+\infty} f^*(t)^pw(t)dt\right)^{\frac{1}{p}}<+\infty\right\},$$
where $f^*$ denotes the non-increasing rearrangement of $f$ (see \cite{BS} for standard notations). Note that if $w=1$ we have $\Lambda^p(w)=L^p$ and if $w(t)=t^{p/q-1}$, with $1\leq q<+\infty$, we obtain $\Lambda^p(w)=L^{q,p}$, where $L^{q,p}$ are the usual Lorentz spaces. In this work we will consider the weighted Lorentz space $\Lambda^p(w)$ such that the weight $w$ satisfies the $B_p$ condition which characterizes the boundedness of the Hardy-Littlewood maximal function on $\Lambda^p(w)$. Indeed, we have $w\in B_p$ for $1\leq p<+\infty$, if there exists $C>0$ such that
$$\int_{r}^{+\infty} \left(\frac{r}{t}\right)^p w(t)dt\leq C \int_{0}^{r}w(t)dt, \mbox{ for all } 0<r<+\infty.$$
and we obtain the inequality $\|\mathscr{M}(f)\|_{\Lambda^p(w)}\leq C \|f\|_{\Lambda^p(w)}$, where $C$ is depending on the quantity
$$[w]_{B_p}=\sup_{r>0}\left\{r^p\, \left(\int_{r}^{+\infty}\frac{w(t)}{t^p}dt \right) \big/ \left( \int_0^r w(t)dt\right)\right\}.$$
For more properties of these weights and the associated classical Lorentz spaces see \cite{Arino}, \cite{Soria} and \cite{Carro}. A generalization of the classical Sobolev inequalities is available in \cite{ChMarcociMarcoci} but the use of rough singular operators seems to be new in the setting of classical Lorentz spaces.\\

In this context, we have the following result.
\begin{Theoreme}\label{Theo_Lorentz}
Over the space $\mathbb{R}^n$ with $n\geq 2$, consider $\Omega$ a function such that $\Omega \in L^1(\mathbb{S}^{n-1})$, $\displaystyle{\int_{\mathbb{S}^{n-1}}\Omega \ d\sigma=0}$ and such that $\Omega\in L^\rho(\mathbb{S}^{n-1})$ with $1<\rho<n$ and consider the operator $T_{\Omega}$ associated to the function $\Omega$ as defined in (\ref{Def_Operator}). \\

\noindent Assume that $f:\mathbb{R}^n\longrightarrow \mathbb{R}$ is a regular function such that $\nabla f\in \dot{\mathcal{M}}^{\mathfrak{p},\mathfrak{q}}(\mathbb{R}^n)$ where $1< \mathfrak{p}\leq \mathfrak{q}<n$ and where $\mathfrak{p}$ is a real parameter such that $1<\overline{\rho}\leq \mathfrak{p}$, with $\overline{\rho}=\frac{\rho n}{\rho n+\rho-n}$. \\

\noindent Consider a weight $w:\mathbb{R}^+\longrightarrow\mathbb{R}^+$  and assume moreover that we have $\nabla f\in \Lambda^{\sigma}(\omega)$ where $\sigma>\mathfrak{p}$. If $\omega\in B_{\frac{r}{\mathfrak{p}}(1-\theta)}$ where $r=\frac{n\sigma}{n-\mathfrak{q}}>1$ and $\theta=\frac{\mathfrak{q}}{n}$, then we have the inequality
\begin{equation}\label{InequalityLorentz}
\|T_\Omega(f)\|_{\Lambda^r(w)}\leq C_\Omega\left\|\nabla f\right\|_{\Lambda^{\sigma}(\omega)}^{1-\theta}\|\nabla f\|_{\dot{\mathcal{M}}^{\mathfrak{p},\mathfrak{q}}}^{\theta}.
\end{equation}
\end{Theoreme}
\noindent {\bf Proof.} Just as in the previous results, we start with the pointwise estimate
$$|T_\Omega(f)(x)|\leq C_\Omega \big(\mathscr{M}(|\nabla f|^{\mathfrak{p}})(x)\big)^{\frac{1}{\mathfrak{p}}-\frac{\mathfrak{q}}{\mathfrak{p}n}}\|\nabla f\|_{\dot{\mathcal{M}}^{\mathfrak{p},\mathfrak{q}}}^{\frac{\mathfrak{q}}{n}}.$$
We will use now the following properties of the non-increasing rearrangement function.
\begin{Lemme}\label{Lem_Lorentz} If $f,g:\mathbb{R}^n\longrightarrow \mathbb{R}$ are two measurable functions, we have
\begin{itemize}
\item[1)] if $|g|\leq |f|$ a.e. then $g^*\leq f^*$,
\item[2)] if $s>0$, then $(|f|^s)^*=(f^*)^s$, in particular we have $\||f^s|\|_{\Lambda^p(\omega)}=\|f\|_{\Lambda^{sp}(\omega)}^s$.
\end{itemize}
\end{Lemme}
\noindent For a proof of this lemma see Proposition 1.4.5 of \cite{Grafakos}. We apply these properties to the previous pointwise estimate to obtain 
\begin{eqnarray*}
|T_\Omega(f)|^*(t)\leq  C_\Omega \left(\mathscr{M}(|\nabla f|^{\mathfrak{p}}\right)^*(t))^{\frac{1}{\mathfrak{p}}-\frac{\mathfrak{q}}{\mathfrak{p}n}}\|\nabla f\|_{\dot{\mathcal{M}}^{\mathfrak{p},\mathfrak{q}}}^{\frac{\mathfrak{q}}{n}}.
\end{eqnarray*}
We multiply now the previous inequality by a weight $w:\mathbb{R}^+\longrightarrow \mathbb{R}^+$ and integrating with respect to the variable $t$, we obtain
$$ \left(\int_{0}^{+\infty}(|T_\Omega(f)|^*(t))^rw(t)dt\right)^{\frac{1}{r}}\leq C_\Omega\|\nabla f\|_{\dot{\mathcal{M}}^{\mathfrak{p},\mathfrak{q}}}^{\frac{\mathfrak{q}}{n}}\left(\int_{0}^{+\infty}\left(\left(\mathscr{M} \left(|\nabla f|^\mathfrak{p}\right)\right)^*(t)\right)^{r(\frac{1}{\mathfrak{p}}-\frac{\mathfrak{q}}{\mathfrak{p}n})}w(t)dt\right)^{\frac{1}{r}}.$$
Since $r(\frac{1}{\mathfrak{p}}-\frac{\mathfrak{q}}{\mathfrak{p}n})>1$ and since $\omega$ belongs to the class $B_{r(\frac{1}{\mathfrak{p}}-\frac{\mathfrak{q}}{\mathfrak{p}n})}$, then we can use the boundedness property of the Hardy-Littlewood maximal function in the spaces $\Lambda^{r(\frac{1}{\mathfrak{p}}-\frac{\mathfrak{q}}{\mathfrak{p}n})}(\omega)$ to obtain 
$$\|T_\Omega(f)\|_{\Lambda^r(w)}\leq C_\Omega\|\nabla f\|_{\dot{\mathcal{M}}^{\mathfrak{p},\mathfrak{q}}}^{\frac{\mathfrak{q}}{n}}\left\||\nabla f|^\mathfrak{p}\right\|_{\Lambda^{r(\frac{1}{\mathfrak{p}}-\frac{\mathfrak{q}}{\mathfrak{p}n})}(\omega)}^{\frac{1}{\mathfrak{p}}-\frac{\mathfrak{q}}{\mathfrak{p}n}}.$$
We can now exploit the last point of the Lemma \ref{Lem_Lorentz} above to obtain
$$\|T_\Omega(f)\|_{\Lambda^r(w)}\leq C_\Omega\|\nabla f\|_{\dot{\mathcal{M}}^{\mathfrak{p},\mathfrak{q}}}^{\frac{\mathfrak{q}}{n}}\|\nabla f\|_{\Lambda^{r(1-\frac{\mathfrak{q}}{n})}(\omega)}^{1-\frac{\mathfrak{q}}{n}},$$
since $r(1-\frac{\mathfrak{q}}{n})=\sigma$ (recall that $r=\frac{n\sigma}{n-\mathfrak{q}}$) and $\theta=\frac{\mathfrak{q}}{n}$, we can deduce the inequality (\ref{InequalityLorentz}) and the Theorem \ref{Theo_Lorentz} is now proven. \hfill $\blacksquare$\\

\noindent {\bf Acknowledgment.} The authors thanks to the referee for careful reading of the paper and for useful suggestions. We also thank Prof. K. Moen for pointing us an error in a previous version of this article. This work was supported by the GDRI ECO-Math. 


\begin{thebibliography}{2}
\bibitem{Arino}
M. \textsc{Ari\~no}, B. \textsc{Muckenhoupt}. \emph{Maximal functions on classical Lorentz spaces and Hardy's inequality with weights for nonincreasing functions}. Trans. Amer. Math. Soc. 320, no. 2, 727--735,  (1990).
\bibitem{BS}
 C. \textsc{Bennett}, R. \textsc{Sharpley.}  \emph{Interpolation of operators}. Pure and Applied Mathematics, Vol 129, Academic Press Inc., Boston, MA, (1988).
\bibitem{Bergh}
J. \textsc{Bergh}, J. \textsc{L\"ofst\"rom}. \emph{Interpolation Spaces}. 
Grundlehren der mathematischen Wissenschaften, 223. Springer Verlag (1976).
\bibitem{Carro}
M. J. \textsc{Carro}, J.A. \textsc{Raposo} and J. \textsc{Soria}. \emph{Recent developments in the theory of Lorentz spaces and weighted inequalities}. Mem. Amer. Math. Soc. 187, no. 877,  (2007).
\bibitem{ChMarcociMarcoci}
D. \textsc{Chamorro}, A. \textsc{Marcoci} and L. \textsc{Marcoci}. \emph{Improved Sobolev inequalities: generalizations to classical Lorentz spaces}. Results in Mathematics, 78:219, (2023). 
\bibitem{Cianchi0}
A. \textsc{Cianchi}. \emph{Strong and weakly inequalities for some classical operators in Orlicz spaces}. J. Lond. Math. Soc. 60(1), 187–202, (1999).
\bibitem{Cianchi1}
A. \textsc{Cianchi}. \emph{Optimal Orlicz-Sobolev embeddings}. Rev. Mat. Iberoamericana 20, no. 2, 427--474 (2004).
\bibitem{Derigoz}
F. \textsc{Deringoz} \emph{et al.} \emph{Generalized fractional maximal and integral operators on Orlicz and generalized Orlicz–Morrey spaces of the third kind}. Positivity 23:727–757, (2019).
\bibitem{Grafakos}
L. \textsc{Grafakos}. \emph{Classical and Modern Fourier Analysis}. Prentice Hall, (2004).
\bibitem{Hoang1}
C. \textsc{Hoang}, K. \textsc{Moen} \& C. \textsc{Pérez}. \emph{New pointwise bounds by Riesz potential type operators}. Journal of Functional Analysis, Volume 289, Issue 9, (2025).
\bibitem{Hoang}
C. \textsc{Hoang}, K. \textsc{Moen} \& C. \textsc{Pérez}. \emph{Pointwise estimates for rough operators with applications to Sobolev inequalities}. Journal d'Analyse Mathématique. Volume 155, 43–74, (2025).
\bibitem{Kinnunen}
J. \textsc{Kinnunen}, J. \textsc{Lehrb\"ack} \& A. \textsc{V\"ah\"akangas}. \emph{Maximal Function Methods for Sobolev Spaces}. Mathematical Surveys and Monographs, Volume 257, American Mathematical Society, (2021).
\bibitem{Li}
K. \textsc{Li}, C. \textsc{Pérez}, I. \textsc{Rivera-Ríos}, and L. \textsc{Roncal}. \emph{Weighted norm inequalities for rough singular integral operators}. J. Geom. Anal. 29, no. 3, 2526–2564, (2019).
\bibitem{Lorentz}
G. G. \textsc{Lorentz}. \emph{Some new functional spaces}, Ann. of Math. (2) 51, 37--55, (1950).
\bibitem{Lorentz1}
G. G. \textsc{Lorentz}. \emph{On the theory of spaces $\Lambda$}, Pacific J. Math. 1, 411--429, (1951).
\bibitem{RaSam}
H. \textsc{Rafeiro}, S. \textsc{Samko}. \emph{Maximal Operator with Rough Kernel in Variable Musielak–Morrey–Orlicz type Spaces, Variable Herz Spaces and Grand Variable Lebesgue Spaces}. Integr. Equ. Oper. Theory 89, 111–124, (2017).
\bibitem{Soria}
J. \textsc{Soria}. \emph{Lorentz spaces of weak-type}. Quart. J. Math. Oxford Ser. (2) 49, no. 193, 93--103, (1998).
\end{thebibliography}
\end{document}